\documentclass[a4paper,pdftex,reqno,10pt]{article}
\usepackage{epsfig,graphicx,times}   
\usepackage{amssymb}
\begin{document}

\title{Counting polyominoes with minimum perimeter}
\author{
Sascha Kurz\thanks{sascha.kurz@uni-bayreuth.de}\\
University of Bayreuth, Department of Mathematics,\\ D-95440 Bayreuth, Germany
}

\maketitle

\noindent
%\begin{abstract}
  \textbf{Abstract.}
  The number of essentially different square polyominoes of order $n$ and minimum perimeter $p(n)$ is enumerated.\\[1mm]
  (In the published version: S. Kurz: Counting polyominoes with minimum perimeter, Ars Combinatoria Vol. 88 (2008), Pages 161-174, 
  there were some problems for the case when the closed walk through the edge-to-edge neighboring squares of the perimeter contains
  a cut vertex. The example from Figure 3 can be treated if the unique square of degree four is counted twice. However, if the unique 
  square of degree four is replaced by a chain of squares of degrees $3$, $2$, $2$, \dots, $2$, $3$, the proofs of Lemma~1 and Lemma~2
  collapse. This does not affect the main result, since for a large number of squares the corresponding graph of polyominoes 
  with minimum perimeter cannot have cut vertices and it can be checked that for a small number of squares the main result is still valid.
  To avoid this inaccuracy, the current version circumvents the auxiliary results of Lemma~1 and Lemma~2.)
%\end{abstract}

\medskip

\noindent
\textbf{Keywords:} polyominoes, enumeration\\
\textbf{MSC:} 05B50$^\star\!$, 05C30\\

\noindent
\rule{\textwidth}{0.3 mm}

  \parskip1.3em
  \section{Introduction}
  Suppose we are given $n$ unit squares. What is the best way to
  arrange them side by side to gain the minimum perimeter $p(n)$? In \cite{HarHar} F. Harary and
  H. Harborth proved that
  $p(n)=2\Big\lceil 2\sqrt{n}\Big\rceil$. They constructed an example
  where the cells grow up cell by cell like spirals for these extremal polyominoes (see Figure 1).
  In general, this is not\\[1mm]
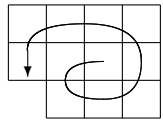
  \begin{center}
  \setlength{\unitlength}{0.5cm}
    \begin{picture}(4,3)
      \put(1,0){\line(1,0){3}}
      \put(0,1){\line(1,0){4}}
      \put(0,2){\line(1,0){4}}
      \put(0,3){\line(1,0){4}}
      \put(0,1){\line(0,1){2}}
      \put(1,0){\line(0,1){3}}
      \put(2,0){\line(0,1){3}}
      \put(3,0){\line(0,1){3}}
      \put(4,0){\line(0,1){3}}
      \qbezier(2.5,1.5)(1.5,1.5)(1.5,1)
      \qbezier(1.5,1)(1.5,0.5)(2.5,0.5)
      \qbezier(2.5,0.5)(3.5,0.5)(3.5,1.5)
      \qbezier(3.5,1.5)(3.5,2.5)(2,2.5)
      \qbezier(2,2.5)(0.5,2.5)(0.5,1.8)
      \put(0.51,1.8){\vector(0,-1){0.7}}
      \put(-2.1,-1.3){Figure 1. Spiral construction.}
    \end{picture}\\[10mm]
  \end{center}
  the only possibility to reach the minimum perimeter.
  Thus the question arises to determine the number $e(n)$ of different square polyominoes
  of order $n$ and with minimum perimeter $p(n)$ where we regard two polyominoes as equal 
  if they can be mapped onto each other by translations, rotations, and reflections.\\
  We will show that these extremal polyominoes can be obtained by
  deleting squares at the corners of rectangular polyominoes with the minimum perimeter $p(n)$
  and with at least $n$ squares. The process of deletion of squares ends if $n$ squares
  remain forming a desired extremal polyomino. This process leads to an enumeration of
  the polyominoes with minimum perimeter $p(n)$.\\[0mm]
  
  \textbf{Theorem 1.}  The number $e(n)$ of polyominoes with $n$ squares and minimum
  perimeter $p(n)$ is given by\\
  \begin{displaymath}
    e(n)=\left\{\begin{array}{l@{\;\quad} c@{\;}l}
    \,\,\,\,\,\,\,\,\,\,\,\,\,\,\,\,\,\,\,1 & if & n=s^2,\\[2mm]
    \sum\limits_{c=0}^{\left\lfloor -\frac{1}{2}+\frac{1}{2}\sqrt{1+4s-4t}\right\rfloor}
     r_{s-c-c^2-t}& if &n=s^2+t,\\
	 && 0< t< s,\\[5mm]
    \,\,\,\,\,\,\,\,\,\,\,\,\,\,\,\,\,\,\,1 & if & n=s^2+s,\\[3mm]
    q_{s+1-t}+\sum\limits_{c=1}^{\left\lfloor
    \sqrt{s+1-t}\right\rfloor}r_{s+1-c^2-t} & if & n=s^2+s+t,\\
	 && 0<t\le s,
    \end{array}\right .
  \end{displaymath}\newpage
  with $s=\lfloor\sqrt{n}\rfloor$, and with $r_k$, $q_k$ being the coefficient of $x^k$ in
  the following generating function $r(x)$ and $q(x)$, respectively. The two generating
  functions
  \begin{displaymath}
    s(x)=1+\sum_{k=1}^\infty x^{k^2}\prod\limits_{j=1}^k\frac{1}{1-x^{2j}}
  \end{displaymath}
  and
  \begin{displaymath}
    a(x)=\prod\limits_{j=1}^\infty\frac{1}{1-x^{j}}
  \end{displaymath}
  
  are used in the definition of\\
  \begin{displaymath}
        r(x)=\frac{1}{4}\left(a(x)^4+3a(x^2)^2\right)
  \end{displaymath}
  and\\[-2mm]
  \begin{center}
    $q(x)=$ $\frac{1}{8}\left(a(x)^4+3a(x^2)^2+2s(x)^2a(x^2)+2a(x^4)\right)\,\,.$
  \end{center}
  \begin{figure}[h]                             % bindet eine Graphik ein
    \centering
    \begin{minipage}[c]{1\textwidth}
      \centering
      Figure 2. $e(n)$ for $n\le 100$.
    \end{minipage}
  \end{figure}
  \noindent
  The behavior of $e(n)$ is illustrated in Figure 2. It has a local
  maximum at $n=s^2+1$ and $n=s^2+s+1$ for $s\ge 1$.
  Then $e(n)$ decreases to $e(n)=1$ at $n=s^2$ and $s=s^2+s$. In the following we give
  lists of the values of $e(n)$ for $n\le 144$ and of the two maximum cases $e(s^2+1)$ and
  $e(s^2+s+1)$ for $s\le 49$,\\[2mm]
  \noindent
  $e(n)=1,$ $1,$ $2,$ $1,$ $1,$ $1,$ $4,$ $2,$ $1,$ $6,$ $1,$ $1,$ $11,$
  $4,$ $2,$ $1,$ $11,$ $6,$ $1,$ $1,$ $28,$ $11,$ $4,$ $2,$ $1,$ $35,$ $11,$
  $6,$ $1,$ $1,$ $65,$ $28,$ $11,$ $4,$ $2,$ $1,$ $73,$ $35,$
  $11,$ $6,$ $1,$ $1,$ $147,$ $65,$ $28,$ $11,$ $4,$ $2,$ $1,$ $182,$
  $73,$ $35,$ $11,$ $6,$ $1,$ $1,$ $321,$ $147,$ $65,$ $28,$ $11,$ $4,$ $2,$ $1,$
  $374,$ $182,$ $73,$ $35,$ $11,$ $6,$ $1,$ $1,$
  $678,$ $321,$ $147,$ $65,$ $28,$ $11,$ $4,$ $2,$ $1,$ $816,$ $374,$ $182,$
  $73,$ $35,$ $11,$ $6,$ $1,$ $1,$ $1382,$ $678,$ $321,$ $147,$ $65,$ $28,$
  $11,$ $4,$ $2,$ $1,$ $1615,$ $816,$ $374,$ $182,$ $73,$ $35,$ $11,$ $6,$ $1,$ $1,$ $2738,$ $1382,$ $678,$
  $321,$ $147,$\newpage $65,$ $28,$ $11,$ $4,$ $2,$ $1,$ $3244,$
  $1615,$ $816,$ $374,$ $182,$ $73,$ $35,$ $11,$ $6,$ $1,$ $1,$ $5289,$ $2738,$ $1382,$ $678,$ $321,$ $147,$
  $65,$ $28,$ $11,$ $4,$ $2,$ $1,$\\[2mm]
  \sloppy
  $e(s^2+1)=1,$ $1,$ $6,$ $11,$ $35,$ $73,$ $182,$ $374,$ $816,$ $1615,$ $3244,$ $6160,$ $11678,$
  $21353,$ $38742,$ $68541,$ $120082,$ $206448,$
  $351386,$ $589237,$ $978626,$ $1605582,$ $2610694,$
  $4201319,$ $6705559,$ $10607058,$ $16652362,$
  $25937765,$ $40122446,$
  $61629301,$ $94066442,$ $142668403,$ $215124896,$ $322514429,$ $480921808,$ $713356789,$\\
  $1052884464,$ $1546475040,$
  $2261006940,$ $3290837242,$ $4769203920,$ $6882855246,$ $9893497078,$ $14165630358,$
  $20206501603,$ $28718344953,$
  $40672085930,$ $57404156326,$ $80751193346,$ \\[2mm]
   $e(s^2+s+1)=2,$ $4,$ $11,$ $28,$ $65,$ $147,$ $321,$ $678,$ $1382,$ $2738,$ $5289,$ $9985,$ $18452,$
  $33455,$ $59616,$ $104556,$ $180690,$
  $308058,$ $518648,$ $863037,$ $1420480,$ $2314170,$ $3734063,$ $5970888,$ $9466452,$
  $14887746,$ $23235296,$
  $36000876,$ $55395893,$ $84680624,$ $128636339,$ $194239572,$ $291620864,$ $435422540,$
  $646713658,$ $955680734,$
  $1405394420,$ $2057063947,$ $2997341230,$ $4348440733,$ $6282115350,$ $9038897722,$ $12954509822,$
  $18496005656,$ $26311093101,$ $37295254695,$ $52682844248,$ $74170401088,$ $104083151128.$
  \section{Proof of Theorem 1}
  The perimeter cannot be a minimum if the polyomino is
  disconnected or if it has holes. For connected polyominoes without holes the property of having
  the minimum perimeter is equivalent to the property of having the maximum number
  of common edges since an edge which does not belong to two squares
  is part of the perimeter. The maximum number of common edges $B(n)$, of a ployomino consisting 
  of $n$ unit squares, is determined in \cite{HarHar} to be $B(n)=2n-\Big\lceil 2\sqrt{n}\Big\rceil$.
  Thus the perimeter of a polyomino consisting of $n$ unit squares and having minimum perimeter
  is given by 
  \renewcommand{\theequation}{$\star$}
  \begin{equation}
    2\Big\lceil 2\sqrt{n}\Big\rceil.    
  \end{equation}

  \textbf{Lemma~1.} For the maximum area $A(H)$ of a polyomino with
  perimeter $H$ we have
  \begin{displaymath}
  A(H)=\left\{
  \begin{array}{l@{\,\,\, if\,\,\,}l} \left(\frac{H}{4}\right)^2 & |H|\equiv
  0\,\,(mod\,4),\\ \left(\frac{H}{4}\right)^2-\frac{1}{4} &
  |H|\equiv 2\,\,(mod\,4).\end{array}\right.\end{displaymath}
  \emph{\textbf{Proof.}}W.l.o.g.\ we can consider a polyomino with minimum 
  perimeter. Because the edges of the perimeter go either in the direction 
  of  the $x$-axis or the direction of the $y$-axis, the integer $H$ has 
  to be an even number. Consider the smallest rectangle surrounding a
  polyomino and denote the side lengths by $a$ and $b$. Using the
  fact that the perimeter $H$ of a polyomino is at least the perimeter of 
  its smallest surrounding rectangle we conclude $H\ge 2a+2b$. The maximum 
  area of the rectangle with given perimeter is obtained if the integers $a$ and
  $b$ are as equal as possible. Thus
  $a=\Big\lceil\frac{H}{4}\Big\rceil$ and
  $b=\Big\lfloor\frac{H}{4}\Big\rfloor$. The product yields
  the asserted formula. \hfill{$\square$}
  
  Now we give a strategy to construct all polyominoes with minimum perimeter. To this 
  end we denote the degree of a square by the number of its edge-to-edge neighbors.
  
  \textbf{Lemma~2.} Each polyomino with minimum perimeter $p(n)$ can be obtained 
  by deleting squares, of degree $2$, of a rectangular polyomino with perimeter 
  $p(n)$ consisting of at least $n$ squares.\\[1.5mm]
  \emph{\textbf{Proof.}}
  Consider a polyomino $P$ with minimum perimeter $H=p(n)$. Denote its smallest 
  surrounding rectangle by $R$. If the perimeter of $R$ is less than $H$ then $P$ 
  does not have the minimum perimeter due to the fact that $m=B(n)$ is increasing.
  Thus $H$ equals the perimeter of $R$ and $P$ can be obtained by deleting squares 
  from a rectangular polyomino with perimeter $p(n)$ and with an area of at least 
  $n$. Only squares of degree~$2$ can be deleted successively if the perimeter does 
  not change.\hfill{$\square$}\\[1mm] 
    
  For the following classes of $n$ with
  $s=\lfloor\sqrt{n}\rfloor$ we now characterize all rectangles being appropriate
  for a deletion process to obtain $P$ with minimum perimeter $p(n)$.\\[5mm]
  \textbf{(i)} $n=s^2$.\\
  From Lemma~1 we know that the unique polyomino with minimum
  perimeter $p(n)$ is indeed the $s\times s$ square.
  
  \noindent
  \textbf{(ii)} $n=s^2+t$, $0<t<s$.\\
  Since \begin{displaymath}
  s^2<n<\left(s+\frac{1}{2}\right)^2=s^2+s+\frac{1}{4}\end{displaymath}
  we conclude from Equation~($\star$) that the perimeter is given by $H=4s+2$. 
  Denote the side lengths of the surrounding rectangle by $a$ and $b$. With
  $2a+2b=H=4s+2$ we let $a=s+1+c$ and $b=s-c$ with an integer
  $c\ge 0$. Since at least $n$ squares are needed for the deletion process
  we have $ab\ge n$,  yielding
  \begin{displaymath} 0\le c\le
  \Big\lfloor-\frac{1}{2}+\frac{1}{2}\sqrt{1+4s-4t}\Big\rfloor.\end{displaymath}
  \textbf{(iii)} $n=s^2+s$.\\
  The $s\times (s+1)$
  rectangle is the unique polyomino with minimum perimeter $p(n)$ due to Lemma~1.
  
  \noindent
  \textbf{(iv)} $n=s^2+s+t$, $0<t\le s$.\\
  Since
  \begin{displaymath}
  \left(s+\frac{1}{2}\right)^2=s^2+s+\frac{1}{4}<n<(s+1)^2=s^2+2s+1
  \end{displaymath}
  we conclude from Equation~($\star$) that the perimeter is given by 
  $H=4s+4$. Again $a$ and $b$ denote the side lengths of the surrounding 
  rectangle and we let $a=s+1+c$ and $b=s+1-c$ with an integer $c\ge 0$. 
  The condition $ab\ge n$ now yields
  \begin{displaymath}0\le c\le
  \Big\lfloor\sqrt{1+s-t}\Big\rfloor.\end{displaymath}
  
  We remark that the deletion process does not change the smallest
  surrounding rectangle since $ab-n < b$, that is the number of deleted
  squares is less than the number of squares of the smallest side of
  this rectangle. (Otherwise the perimeter would decrease.)
  
  So far we have described those rectangles from which squares of degree 2 are removed.
  Now we examine the process of deleting squares from
  a rectangular polyomino. Squares
  of degree 2 can only be located in the corners of the polyomino.
  What shape has the set of deleted squares at any corner?
  There is a maximum square\\[3mm]
  \setlength{\unitlength}{0.5cm}
  \begin{picture}(21,4)
    % First example
    \put(0,0){\line(0,1){4}}
    \put(1,0){\line(0,1){4}}
    \put(2,0){\line(0,1){4}}
    \put(3,1){\line(0,1){3}}
    \put(4,1){\line(0,1){3}}
    \put(0,0){\line(1,0){2}}
    \put(0,1){\line(1,0){4}}
    \put(0,2){\line(1,0){6}}
    \put(0,3){\line(1,0){8}}
    \put(0,4){\line(1,0){8}}
    \put(8,4){\line(0,-1){1}}
    \put(7,4){\line(0,-1){1}}
    \put(6,4){\line(0,-1){2}}
    \put(5,4){\line(0,-1){2}}
    \thicklines
    \put(0,1){\line(1,0){3}}
    \put(0,1){\line(0,1){3}}
    \put(3,4){\line(-1,0){3}}
    \put(3,4){\line(0,-1){3}}
    %\put(0,4){\line(1,-1){3}}
    % Second example
    \thinlines
    \put(10,0){\line(1,0){1}}
    \put(10,3){\line(1,0){5}}
    \put(12,4){\line(1,0){3}}
    \put(12,2){\line(1,0){2}}
    \put(10,0){\line(1,0){1}}
    \put(10,0){\line(0,1){2}}
    \put(11,0){\line(0,1){4}}
    \put(13,4){\line(0,-1){2}}
    \put(14,4){\line(0,-1){2}}
    \put(15,4){\line(0,-1){1}}
    \put(10,1){\line(1,0){1}}
    \thicklines
    \put(10,2){\line(1,0){2}}
    \put(10,2){\line(0,1){2}}
    \put(12,4){\line(-1,0){2}}
    \put(12,4){\line(0,-1){2}}
    %\put(10,4){\line(1,-1){2}}
    % Third example
    \thinlines
    \put(17,0){\line(1,0){1}}
    \put(17,1){\line(1,0){2}}
    \put(17,0){\line(0,1){2}}
    \put(18,0){\line(0,1){4}}
    \put(19,1){\line(0,1){1}}
    \put(20,2){\line(0,1){2}}
    \put(21,3){\line(0,1){1}}
    \put(21,4){\line(-1,0){2}}
    \put(21,3){\line(-1,0){4}}
    \put(20,2){\line(-1,0){1}}
    \thicklines
    \put(17,2){\line(1,0){2}}
    \put(17,2){\line(0,1){2}}
    \put(19,4){\line(-1,0){2}}
    \put(19,4){\line(0,-1){2}}
    %\put(17,4){\line(1,-1){2}}
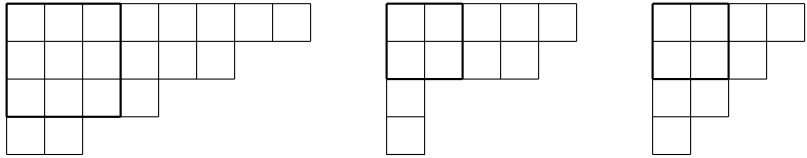
    \put(2,-1.3){Figure 4. Shape of the deleted squares at the corners.}
  \end{picture}\\[10mm]
  of squares at the corner, the so called \lq\lq Durfee square", together with squares in rows and columns
  of decreasing length from
  outside to the interior part of the polyomino. To count the different possibilities
  of the sets of deleted squares with respect to the number of the deleted squares
  we use the concept of a generating function
  $f(x)=\sum_{i=0}^{\infty}f_ix^i$. Here the coefficient $f_i$ gives the number of
  different ways to use $i$ squares. Since the rows and columns
  are ordered by their lengths they form Ferrer's diagrams with generating function
  $\prod_{j=1}^\infty\frac{1}{1-x^j}$ each \cite{Partition}. So the generating function for the sets of 
  deleted squares in a single corner is given by
  \begin{displaymath}
    a(x)=\prod_{j=1}^\infty\frac{1}{1-x^j}\,. 
  \end{displaymath}
  Later we will also need the generating function $s(x)$ for the sets of deleted squares 
  being symmetric with respect to the diagonal of the corner square. Since such a symmetric set of deleted
  squares consists of a square of $k^2$ squares and the two mirror images of a 
  Ferrer's diagrams with height or width at most $k$ we get
  \begin{displaymath}
    s(x)=1+\sum_{k=1}^\infty x^{k^2}\prod_{j=1}^k\frac{1}{1-x^{2j}}\,.
  \end{displaymath}
  We now consider the whole rectangle. Because of different sets
  of symmetry axes we distinguish between squares and rectangles.
  We define generating functions $q(x)$ and $r(x)$ so that the coefficient of $x^k$
  in $q(x)$ and $r(x)$ is the number of ways to remove $k$
  squares from all four corners of a square or a rectangle, respectively. We mention that
  the coefficient of
  $x^k$ gives the desired number only if $k$ is smaller than the small
  side of the rectangle.\\[4mm]
  Since we want to count polyominoes with minimum perimeter up to translation,
  rotation, and reflection, we have to factor out these symmetries. Here the
  general tool is the lemma of Cauchy-Frobenius, see e.g. \cite{CycleInd}. We remark 
  that we do not have to consider translations because we describe the polyominoes 
  without coordinates.
  
  \textbf{Lemma (Cauchy-Frobenius, weighted form).} Given a group action of a finite group $G$ on a set $S$ 
  and a map $w:S\longrightarrow R$ from $S$ into a commutative ring $R$ containing 
  $\mathbb{Q}$ as a subring. If $w$ is constant on the orbits of $G$ on $S$, then we have, for any transversal 
  $\mathcal{T}$ of the orbits:
  $$
  	\sum_{t\in\mathcal{T}}w(t)=\frac{1}{|G|}\sum_{g\in G}\sum_{s\in S_g}w(s)
  $$
  where $S_g$ denotes the elements of $S$ being fixed by $g$, i.e. 
  $$
    S_g=\lbrace s\in S|s=gs\rbrace \,.
  $$

  For $G$ we take the symmetry group of a square or a rectangle, respectively, for $S$ we take the sets 
  of deleted squares on all 4 corners, and for the weight $w(s)$ we take $x^k$, where $k$ is the number of 
  squares in $s$.
  Here we will only describe in detail the application of this lemma for a determination of $q(x)$. We label 
  the 4 corners of the square by $1$, $2$, $3$, and $4$, see Figure~5. In Table~1 we list the 8\\
  
  \begin{center}
  \setlength{\unitlength}{0.5cm}
  \begin{picture}(12.3,5.4)
    \put(0,0){$\mathbf{4}$}
	 \put(4.7,0){$\mathbf{3}$}
	 \put(0,4.8){$\mathbf{1}$}
	 \put(4.6,4.8){$\mathbf{2}$}
    \put(0.6,0.7){\line(1,0){4}}
    \put(0.6,0.7){\line(0,1){4}}
    \put(4.6,4.7){\line(-1,0){4}}
    \put(4.6,4.7){\line(0,-1){4}}
    \put(0.6,3.7){\line(1,0){2}}
    \put(1.6,3.7){\line(0,1){1}}
    \put(2.6,3.7){\line(0,1){1}}

	 \put(7,0){$\mathbf{4}$}
	 \put(11.7,0){$\mathbf{3}$}
	 \put(7,4.8){$\mathbf{1}$}
	 \put(11.6,4.8){$\mathbf{2}$} 
    \put(7.6,0.7){\line(1,0){4}}
    \put(7.6,0.7){\line(0,1){4}}
    \put(11.6,4.7){\line(-1,0){4}}
    \put(11.6,4.7){\line(0,-1){4}}
    \put(8.6,4.7){\line(0,-1){2}}
    \put(7.6,2.7){\line(1,0){1}}
    \put(7.6,3.7){\line(1,0){1}}
  \end{picture}\\[2mm]
  Figure 5. 
  \end{center}
  
  permutations $g$ of 
  the symmetry group of a square, the dihedral group on 4 points, together with the 
  corresponding generating functions for the sets $S_g$ being fixed by $g$.
  
  \begin{center}
    \begin{tabular}{cc}
    $(1)(2)(3)(4)$ & $a(x)^4$ \\
    $(1,2,3,4)$ & $a(x^4)$ \\
    $(1,3)(2,4)$&$a(x^2)^2$\\
    $(1,4,3,2)$&$a(x^4)$  \\
    $(1,2)(3,4)$&$a(x^2)^2$  \\
    $(1,4)(2,3)$&$a(x^2)^2$  \\
    $(1,3)(2)(4)$&$s(x)^2a(x^2)$  \\
    $(1)(2,4)(3)$&$s(x)^2a(x^2)$  \\
    \end{tabular}\\[2mm]
    Table 1. Permutations of the symmetry group of a square together with the
    corresponding generating functions of $S_g$.\\[15mm]
  \end{center}
  
  The generating function of the set of deleted squares on a corner is $a(x)$. If we consider 
  the configurations being fixed by the identity element $(1)(2)(3)(4)$ we see that the sets of deleted squares 
  at the 4 corners are independent and so $|S_{(1)(2)(3)(4)}|=a(x)^4$. In the case when $g=(1,2,3,4)$ the sets of 
  deleted squares have to be the same for all 4 corners and we have $|S_{(1,2,3,4)}|=a(x^4)$. For the double transposition 
  $(1,2)(3,4)$ the sets of deleted squares at corners $1$ and $2$, and the sets of deleted squares at corners 
  $3$ and $4$ have to be equal. Because the sets of deleted squares at corner points $1$ and $3$ are independent we get 
  $|S_{(1,2)(3,4)}|=a(x^2)^2$. Next we consider $g=(1)(2,4)(3)$. The sets of deleted squares at corners $2$ and 
  $4$ have to be equal. If we apply $g$ on the polyomino of the left hand side of Figure 5 we receive the polyomino on the
  right hand side and we see that in general the sets of deleted squares at corners $1$ and $3$ have to be symmetric. 
  Thus $|S_{(1)(2,4)(3)}|=s(x)^2a(x^2)$. The other cases are left to the reader. Summing up and a division by 8 yields   
   \begin{center}
    $q(x)=$ $\frac{1}{8}\left(a(x)^4+3a(x^2)^2+2s(x)^2a(x^2)+2a(x^4)\right)\,\,.$
  \end{center}
  \noindent
  For the symmetry group of a rectangle we analogously obtain
  \begin{displaymath}
    r(x)=\frac{1}{4}\left(a(x)^4+3a(x^2)^2\right).
  \end{displaymath}
  
  With the preceding characterization of rectangles being appropriate for a deletion process and the 
  formulas for $a(x)$, $s(x)$, $q(x)$, and $r(x)$ we have the proof of Theorem~1 at hand.\\[8mm]
  We would like to close with the first entries of a complete list of polyominoes with minimum perimeter $p(n)$, 
  see Figure~6.\\[-2mm]
  
  \section{Acknowledgments}
  I thank Ralf Gugisch, Heiko Harborth, Adalbert Kerber, and Axel Kohnert for many fruitful discussions 
  during the preparation of this article.
  
  \begin{center}
  \setlength{\unitlength}{0.5cm}
  \begin{picture}(22,13)
    % 4 x 3 Polyomino minus left upper corner and the corner below
    \put(8,0){\line(1,0){4}}
    \put(8,0){\line(0,1){3}}
    \put(8,3){\line(1,0){3}}
    \put(11,1){\line(0,1){2}}
    \put(11,1){\line(1,0){1}}
    \put(12,0){\line(0,1){1}}
    % 4 x 3 Polyomino minus left upper corner and the corner below
    % 3 x 4 Polyomino minus left upper and minus right upper corner
    \put(4,0){\line(1,0){3}}
    \put(4,0){\line(0,1){3}}
    \put(4,3){\line(1,0){1}}
    \put(5,3){\line(0,1){1}}
    \put(5,4){\line(1,0){1}}
    \put(6,4){\line(0,-1){1}}
    \put(6,3){\line(1,0){1}}
    \put(7,0){\line(0,1){3}}
    % 3 x 4 Polyomino minus left upper and minus right upper corner
    % 3 x 4 Polyomino minus left upper corner and the corner below
    \put(0,0){\line(1,0){3}}
    \put(0,0){\line(0,1){4}}
    \put(0,4){\line(1,0){2}}
    \put(2,4){\line(0,-1){2}}
    \put(2,2){\line(1,0){1}}
    \put(3,0){\line(0,1){2}}
    % 3 x 4 Polyomino minus left upper corner and the corner below
    % 4 x 3 Polyomino minus left lower and minus right lower corner
    \put(15.5,6){\line(1,0){4}}
    \put(15.5,4){\line(0,1){2}}
    \put(15.5,4){\line(1,0){1}}
    \put(16.5,3){\line(0,1){1}}
    \put(16.5,3){\line(1,0){2}}
    \put(18.5,3){\line(0,1){1}}
    \put(18.5,4){\line(1,0){1}}
    \put(19.5,4){\line(0,1){2}}
    % 4 x 3 Polyomino minus left lower and minus right lower corner
    % 2 x 5 Polyomino
    \put(8,4){\line(1,0){5}}
    \put(8,4){\line(0,1){2}}
    \put(13,4){\line(0,1){2}}
    \put(8,6){\line(1,0){5}}
    % 2 x 5 Polyomino
    % 4 x 3 Polyomino minus left lower and minus right upper corner
    \put(14,0){\line(1,0){3}}
    \put(14,0){\line(0,1){1}}
    \put(14,1){\line(-1,0){1}}
    \put(13,1){\line(0,1){2}}
    \put(13,3){\line(1,0){3}}
    \put(16,3){\line(0,-1){1}}
    \put(16,2){\line(1,0){1}}
    \put(17,2){\line(0,-1){2}}
    % 4 x 3 Polyomino minus left lower and minus right upper corner
    % 1 x 1 Polyomino
    \put(0,12){\line(1,0){1}}
    \put(0,12){\line(0,1){1}}
    \put(1,13){\line(-1,0){1}}
    \put(1,13){\line(0,-1){1}}
    % 1 x 1 Polyomino
    % 2 x 3 Polyomino
    \put(0,9){\line(1,0){3}}
    \put(0,9){\line(0,1){2}}
    \put(3,11){\line(-1,0){3}}
    \put(3,11){\line(0,-1){2}}
    % 2 x 3 Polyomino
    % 3 x 3 Polyomino minus right upper corner
    \put(0,5){\line(1,0){3}}
    \put(0,5){\line(0,1){3}}
    \put(0,8){\line(1,0){2}}
    \put(3,5){\line(0,1){2}}
    \put(2,8){\line(0,-1){1}}
    \put(2,7){\line(1,0){1}}
    % 3 x 3 Polyomino minus right upper corner
    % 2 x 1 Polyomino
    \put(2,12){\line(1,0){2}}
    \put(2,12){\line(0,1){1}}
    \put(4,13){\line(-1,0){2}}
    \put(4,13){\line(0,-1){1}}
    % 2 x 1 Polyomino
    % 3 x 1 Polyomino
    \put(5,12){\line(1,0){3}}
    \put(5,12){\line(0,1){1}}
    \put(8,13){\line(-1,0){3}}
    \put(8,13){\line(0,-1){1}}
    % 3 x 1 Polyomino
    % 2 x 2 Polyomino minus right lower corner
    \put(9,11){\line(1,0){1}}
    \put(9,11){\line(0,1){2}}
    \put(9,13){\line(1,0){2}}
    \put(11,13){\line(0,-1){1}}
    \put(11,12){\line(-1,0){1}}
    \put(10,12){\line(0,-1){1}}
    % 2 x 2 Polyomino minus right lower corner
    % 2 x 2 Polyomino
    \put(12,11){\line(1,0){2}}
    \put(12,11){\line(0,1){2}}
    \put(14,13){\line(-1,0){2}}
    \put(14,13){\line(0,-1){2}}
    % 2 x 2 Polyomino
    % 3 x 2 Polyomino minus right lower corner
    \put(15,11){\line(1,0){2}}
    \put(15,11){\line(0,1){2}}
    \put(15,13){\line(1,0){3}}
    \put(18,13){\line(0,-1){1}}
    \put(18,12){\line(-1,0){1}}
    \put(17,12){\line(0,-1){1}}
    % 3 x 2 Polyomino minus right lower corner
    % 3 x 3 Polyomino
    \put(4,5){\line(1,0){3}}
    \put(4,5){\line(0,1){3}}
    \put(7,8){\line(-1,0){3}}
    \put(7,8){\line(0,-1){3}}
    % 3 x 3 Polyomino
    % 4 x 2 Polyomino minus right lower corner
    \put(4,9){\line(1,0){3}}
    \put(4,9){\line(0,1){2}}
    \put(4,11){\line(1,0){4}}
    \put(8,11){\line(0,-1){1}}
    \put(8,10){\line(-1,0){1}}
    \put(7,10){\line(0,-1){1}}
    % 4 x 2 Polyomino minus right lower corner
    % 3 x 3 Polyomino minus left upper and minus right lower corner
    \put(8,7){\line(1,0){2}}
    \put(8,7){\line(0,1){2}}
    \put(10,7){\line(0,1){1}}
    \put(10,8){\line(1,0){1}}
    \put(11,8){\line(0,1){2}}
    \put(11,10){\line(-1,0){2}}
    \put(9,10){\line(0,-1){1}}
    \put(9,9){\line(-1,0){1}}
    % 3 x 3 Polyomino minus left upper and minus right lower corner
    % 3 x 3 Polyomino minus right upper and minus right central corner
    \put(12,7){\line(1,0){3}}
    \put(12,7){\line(0,1){3}}
    \put(15,7){\line(0,1){1}}
    \put(12,10){\line(1,0){2}}
    \put(14,10){\line(0,-1){2}}
    \put(14,8){\line(1,0){1}}
    % 3 x 3 Polyomino minus right upper and minus right central corner
    % 3 x 3 Polyomino minus right upper and minus right lower corner
    \put(16,7){\line(1,0){2}}
    \put(16,7){\line(0,1){3}}
    \put(18,7){\line(0,1){1}}
    \put(18,8){\line(1,0){1}}
    \put(19,8){\line(0,1){1}}
    \put(19,9){\line(-1,0){1}}
    \put(18,9){\line(0,1){1}}
    \put(18,10){\line(-1,0){2}}
    % 3 x 3 Polyomino minus right upper and minus right lower corner
    % 2 x 4 Polyomino
    \put(20,11){\line(1,0){2}}
    \put(20,11){\line(0,-1){4}}
    \put(20,7){\line(1,0){2}}
    \put(22,7){\line(0,1){4}}
    % 2 x 4 Polyomino
    % 3 x 4 Polyomino minus left upper corner
    \put(20,4){\line(1,0){2}}
    \put(20,4){\line(0,-1){1}}
    \put(22,4){\line(0,-1){4}}
    \put(20,3){\line(-1,0){1}}
    \put(19,3){\line(0,-1){3}}
    \put(19,0){\line(1,0){3}}
    % 3 x 4 Polyomino minus left upper corner
  \end{picture}\\[2mm]
  \
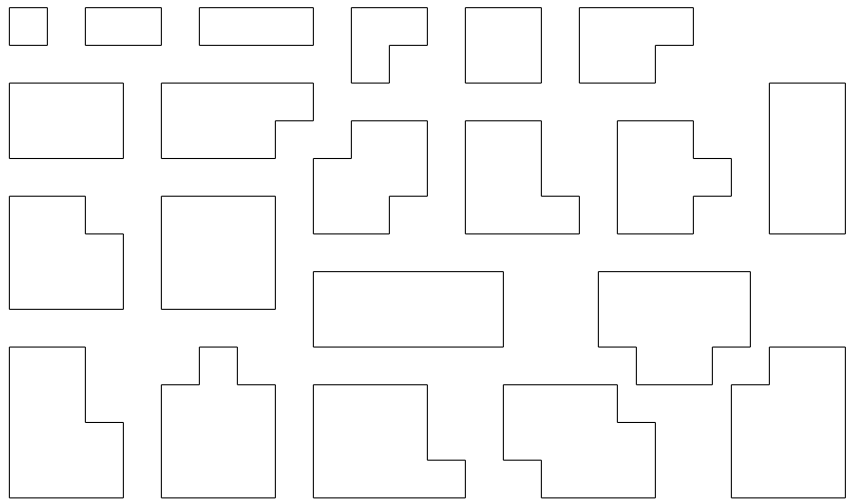
  Figure 6. Polyominoes with minimum perimeter $p(n)$ for $n\le 11$.
  \end{center}

  \nocite{Cube}
  \bibliography{extremal_polyominoes_paper1}
  \bibdata{extremal_polyominoes_paper1}
  \bibliographystyle{amsplain}
\end{document}